\documentstyle[12pt]{article}
 \newcommand{\sect}[1]{\setcounter{equation}{0}\section{#1}}

 \font\tengoth=eufm10 \font\sevengoth=eufm7 \font\fivegoth=eufm5
 \newfam\gothfam
 \textfont\gothfam=\tengoth \scriptfont\gothfam=\sevengoth
   \scriptscriptfont\gothfam=\fivegoth
 \font\frak=eufm10 scaled\magstep1
 \newcommand{\fra}[1]{\mbox{\frak #1}}

 \def\be{\begin{equation}}
 \def\ee{\end{equation}}
 \def\bea{\begin{eqnarray}}
 \def\eea{\end{eqnarray}}
 
 \def\g{{\fra g}}

 \def\unal{\mbox{$U({\bf L})$ }}
 \def\unalj{\mbox{$U_j({\bf L})$ }}
 \def\unfial{\mbox{$U_{j}({\bf L})$ }}
 \def\unfal{\mbox{$U_{E}({\bf L})$ }}

 \def\delfih{\mbox{$\Delta_{j}(H)$}}
 \def\delfia{\mbox{$\Delta_{j}(A)$}}
 \def\delfib{\mbox{$\Delta_{j}(B)$}}
 \def\delfie{\mbox{$\Delta_{j}(E)$}}

 \def\delfh{\mbox{$\Delta_{E}(H)$}}
 \def\delfa{\mbox{$\Delta_{E}(A)$}}
 \def\delfb{\mbox{$\Delta_{E}(B)$}}
 \def\delfe{\mbox{$\Delta_{E}(E)$}}

 \def\delfcc{\mbox{$\Delta_{P'}$}}
 \def\delfcha{\mbox{$\Delta_{P}(H)$}}
 \def\delfcaa{\mbox{$\Delta_{P}(A)$} }
 \def\delfcba{\mbox{$\Delta_{P}(B)$}}
 \def\delfcea{\mbox{$\Delta_{P}(E)$}}

 \def\exs{\mbox{$e^{\sigma}$}}
 \def\exms{\mbox{$e^{-\sigma}$}}
 \def\ex2s{\mbox{$e^{2\sigma}$}}
 \def\exm2s{\mbox{$e^{-2\sigma}$}}
 \def\exas{\mbox{$e^{\alpha \sigma}$}}
 \def\exbes{\mbox{$e^{\beta \sigma}$}}
 \def\exdes{\mbox{$e^{\delta \sigma}$}}

 \def\exmbes{\mbox{$e^{-\beta \sigma}$}}
 \def\exmdes{\mbox{$e^{-\delta \sigma}$}}
  
\begin{document}

\ \vskip 1cm

\begin{center}
 
{\LARGE{\bf{Peripheric Extended Twists}}}
 
\vskip1cm
 
{\sc Vladimir Lyakhovsky \footnote{On leave of absence from Theoretical Department,
  Sankt-Petersburg State University, 198904, St. Petersburg, Russia.

This work has been partially
supported by DGES of the  Ministerio de  Educaci\'on y Cultura of Espa\~na under Project
PB95-0719, the Junta de Castilla y Le\'on (Espa\~na) and the Russian 
Foundation for Fundamental Research under grant  97-01-01152.
} and Mariano A. del Olmo} 
\vskip0.5cm
 
{\it   Departamento de  F\'{\i}sica Te\'orica,
  Facultad de Ciencias}\\
 {\it Universidad de Valladolid,  E-47011, Valladolid, Spain.}
\vskip0.25cm
 
 {e.mail: lyakhovsky@fta.uva.es ;
olmo@fta.uva.es}
 \end{center}
 
\vskip1cm
\centerline{\today}
\vskip1.75cm
 
\begin{abstract}
The properties of the set $\widehat{\cal L}$ of extended jordanian twists 
are studied.
It is shown that the boundaries of $\widehat{\cal L}$ contain twists  
whose characteristics differ considerably from those of internal points.
The extension multipliers of these ``peripheric"  twists are factorizable. 
This leads to simplifications in the twisted algebra relations and helps to find
the explicit form for coproducts. The peripheric twisted algebra $U(sl(4))$
is obtained to illustrate the construction. It is shown that
the corresponding deformation $U_{P}(sl(4))$ cannot be connected 
with the Drinfeld--Jimbo one by a smooth limit procedure. All the 
carrier algebras for the extended and the peripheric extended twists are 
proved to be Frobenius.     
\end{abstract}

\newpage

\sect{Introduction}

Any Lie bialgebra has a quantum deformation \cite{ETI}.
Though there are not so many cases when it can be written in the
global form. In this context the explicit knowledge of the universal
${\cal R}$--matrix is of crucial
importance. It provides the possibility  to build  the
$R$--matrices in any representation and to use the advantages of
the FRT formalism  \cite{FAD}. This is why the triangular Hopf algebras and
twists (they preserve the triangularity \cite{D2,D3}) play such 
an important role in quantum
group theory and applications \cite{KUL1,VLA,VAL}. Despite these
facts very few types of  twists were written explicitly in a
closed form. The well known example is the jordanian twist (JT) of
$sl(2)$ or, more exactly, of its Borel
subalgebra $B(2)$ ($\{H,E|[H,E]=2E\}$)
with $r=H\otimes E-E\otimes H=H\wedge E$ \cite{DRIN} where the 
triangular $R$%
-matrix ${\cal R}=({\cal F}_j)_{21}{\cal F}_j^{-1}$  is defined by the twisting
element \cite{OGIEV,GER}
\begin{equation}
\label{og-twist}
{\cal F}_j=\exp \{\frac 12H\otimes \ln (1+2\xi E)\}.
\end{equation}
In \cite{KLM} it was shown that there exist different extensions (ET's) of this twist.
In particular the ET deformation for ${\cal U}
(sl(N))$ was constructed with the explicit expressions of deformed compositions. 
Using
the notion of factorizable twist \cite{RSTS} the element $%
{\cal F}_E\in {\cal U} (sl(N))^{\otimes 2}$,
\begin{equation}
\label{twist-sl(N)}
{\cal F}_E=\exp \{2\xi \sum_{i=2}^{N-1}E_{1i}\otimes
E_{iN}e^{-\sigma }\}\exp \{H\otimes \sigma \},
\end{equation}
was proved to satisfy the twist equation,
where $E=E_{1N}$, $H=E_{11}-E_{NN}$ and $\sigma =\frac 12\ln (1+2\xi E)$. For 
simplicity of compositions the algebra $sl(N)$ is presented above in the 
standard $gl(N)$ basis: $\{ E_{ij} \} _{i,j = 1, \dots N} $.

The connection of the Drinfeld--Jimbo (DJ) deformation \cite{DRIN,JIMB} with the
jordanian deformation was already pointed out in \cite{GER}. The similarity
transformation of the classical matrix 
$$
r_{DJ} = \sum_{i=1}^{{\rm rank}(g)}t_{ij}H_i \otimes H_j + \sum_{\alpha \in \Delta_+}
E_{\alpha} \otimes E_{- \alpha}
$$ 
performed by the operator $\exp (\xi
{\rm ad}E_{1N})$ (with the highest root generator $E_{1N})$ turns $r_{DJ}$
into the sum $r_{DJ}+\xi r_j$ \cite{GER} where 
\begin{eqnarray}
r_j=-\xi \left( H_{1N}\wedge E_{1N}+2\sum_{k=2}^{N-1}E_{1k}\wedge
E_{kN}\right) .
\end{eqnarray}
Hence $r_j$ is also a classical $r$--matrix  and  defines the corresponding
deformation. A contraction of the quantum
Manin plane $xy=qyx$ of ${\cal U}_q(sl(2))$
with the mentioned above similarity transformation in the fundamental
representation $M=1+\theta \rho (E_{12})$, $\theta =\xi
(1-q)^{-1}$ results
in the jordanian plane $x^{\prime }y^{\prime }
=y^{\prime }x^{\prime }+\xi {y^{\prime }}^2$ of ${\cal U}_{j}(sl(2))$
\cite{OGIEV}. Thus, the jordanian and the extended 
jordanian twisted algebras (with the carrier subalgebra correlated with the
standard dual $\g^*_{DJ}$) can be treated as a limit case for the 
parameterized set of Drinfeld--Jimbo quantizations.

In this paper we study the family of carrier algebras 
(the term is considered to appear first in \cite{GER}) of the type  ${\bf
L}$,  that is the three-parametric set 
${\cal L}=\{{\bf L}(\alpha,\beta,\gamma,\delta)_{\alpha + \beta = \delta}\}$ 
and the  properties of the corresponding sets $\widehat{\cal L}$ of twists 
when the parameters tend to
its  limit values (Section 3). We show that there are two cases ($\alpha
\rightarrow 0 $ and $\beta \rightarrow 0 $) when the twists survive and 
remain nontrivial. We call these twists  peripheric extended twists
(PE twists or PET's), they form the boundary subsets of the variety
$\widehat{\cal L}$. 

The properties of the peripheric algebras differ considerably from those of the 
internal points of ${\cal L}$. The same is true for the properties of 
PE twists. Contrary to the general situation the extension factors of PE 
twists are the solutions of the factorized twist equations (see Section 2).
In Section 4 we show how ${\bf L}(0,\beta,\gamma,\beta)$ or
${\bf L}(\alpha,0,\gamma,\alpha)$ can be injected in the simple Lie algebras 
and illustrate all the results for the case $sl(4) \supset 
{\bf L}(-1,0,1,-1)$. The  deformed coproducts thus obtained for
$U_{{\cal F}_P}(sl(4))$ are much simpler than in the case of general ET's
and the complete list of them for the generators of $sl(4)$ is presented.
The other significant fact is that the PE twists cannot be connected with the
DJ deformations by any kind of smooth ``contraction" (Section 5). The
solutions of the classical Yang--Baxter equation corresponding to ET's and 
PET's can be easily related with the classification given by Stolin
\cite{STO}.  The internal points of the variety ${\cal L}$ are Frobenius
algebras. On the boundary
only the above mentioned subsets $\{ {\bf L}(0,\beta,\gamma,\beta) \}$ and 
$\{{\bf L}(\alpha,0,\gamma,\alpha) \}$ are formed by Frobenius algebras. 
The paper is concluded by the discussion of relations between 
Drinfeld-Jimbo, extended twist and peripheric extended twist
deformations. 


\sect{Basic definitions }

In this section we remind briefly the basic
notions connected with the twisting procedure.

A Hopf algebra ${\cal A}(m,\Delta ,\epsilon,S)$ with
multiplication $m\colon {\cal A}\otimes {\cal A}\to {\cal A}$,
coproduct $\Delta \colon {\cal A}\to {\cal A}%
\otimes {\cal A}$, counit $\epsilon \colon {\cal A}\to C$,
and antipode $S : {\cal A}\to {\cal A}$ 
can be transformed \cite{D2} with an invertible (twisting) 
element ${\cal F}\in {\cal A} 
\otimes {\cal A}$, ${\cal F}=\sum f_i^{(1)}\otimes f_i^{(2)}$, 
into a twisted
one ${\cal A}_{\cal F}(m,\Delta _{\cal F},\epsilon ,S_{\cal F})$.
This Hopf algebra ${\cal A}_{\cal F}$ has the
same multiplication and counit  but the twisted coproduct and antipode:
\begin{equation}
\label{def-t}
\Delta _{\cal F}(a)={\cal F}\Delta (a){\cal F}^{-1},\qquad S_{\cal
F}(a)=vS(a)v^{-1},
\end{equation}
with
$$
v=\sum f_i^{(1)}S(f_i^{(2)}), \qquad a\in {\cal A}.
$$
The twisting element has to satisfy the equations
\begin{eqnarray}
\label{def-n}
(\epsilon \otimes  id)({\cal F}) = (id \otimes  \epsilon)({\cal F})=1,\\[0.2cm]
\label{gentwist}
{\cal F}_{12}(\Delta \otimes  id)({\cal F}) =
{\cal F}_{23}(id \otimes  \Delta)({\cal F});
\label{TE}
\end{eqnarray}
the first one is just a normalization condition and
follows from the second relation modulo a non-zero scalar factor.

If ${\cal A}$ is a Hopf subalgebra of ${\cal B}$ the twisting 
element ${\cal F}$
satisfying (\ref{def-t})--(\ref{gentwist}) induces the twist
deformation  ${\cal B}_{\cal F}$ of  ${\cal B}$. In this case one can
put $a \in  {\cal B}$ in all the formulas (\ref{def-t}). This will
completely define the Hopf algebra ${\cal B}_{\cal F}$. Let ${\cal A}$ and
${\cal B}$ be the universal enveloping algebras: ${\cal A} = U({\fra l}) \subset
{\cal B}= U(\g)$ with ${\fra l} \subset \g$. If $U({\fra l})$ is the minimal
subalgebra on which ${\cal F}$ is completely defined as ${\cal F} \in U({\fra l})
\otimes U({\fra l})$ then ${\fra l}$ is called the carrier algebra for ${\cal F}$
\cite{GER}. 

The composition of appropriate twists can be defined as
${\cal F} = {\cal F}_2 {\cal F}_1$. Here the element ${\cal F}_1$ has to
satisfy the twist equation with the coproduct of the original Hopf algebra,
while ${\cal F}_2$ must be its solution for $\Delta_{{\cal F}_1}$ of the
one twisted by ${\cal F}_1$.
In particular, if ${\cal F}$ is a solution to the twist equation (\ref{TE})
then ${\cal F}^{-1}$ satisfies this equation with $\Delta$ substituted by 
$\Delta_{\cal F}$.

If the initial Hopf algebra ${\cal A}$ is quasitriangular with the 
universal element ${\cal R}$ then so is the twisted one 
${\cal A}_{\cal F}(m,\Delta _{\cal F},\epsilon,S_{\cal F},{\cal R}_{\cal F})$
whose universal element is related to the initial ${\cal R}$ by a transformation
\begin{eqnarray}\label{Rt}
{\cal R}_{\cal F}=  {\cal F}_{21} \,{\cal R} \,{\cal F}^{-1}.
\end{eqnarray}

Most of the explicitly known twisting elements have the factorization
property with respect to comultiplication
$$
(\Delta \otimes id)({\cal F})={\cal F}_{23}{\cal F}_{13}\qquad \mbox{or}
\qquad
(\Delta \otimes id)({\cal F})={\cal F}_{13}{\cal F}_{23}\,,
$$
and
$$
(id \otimes \Delta)({\cal F})={\cal F}_{12}{\cal F}_{13}\qquad \mbox{or}
\qquad
(id \otimes \Delta)({\cal F})={\cal F}_{13}{\cal F}_{12}\,.
$$
To guarantee the validity of the twist equation, these identities are to be combined with 
the additional requirement
${\cal F}_{12}{\cal F}_{23}={\cal F}_{23}{\cal F}_{12}$
or the Yang--Baxter equation on ${\cal F}$ \cite{RSTS}.

An important subclass of factorizable twists consists of elements
satisfying the equations
\begin{eqnarray} \label{f-twist1}
(\Delta \otimes id)({\cal F})={\cal F}_{13}{\cal F}_{23}\,,   
\\ [0,2cm] \label{f-twist2}
(id\otimes \Delta _{\cal F})({\cal F})={\cal F}_{12}{\cal F}_{13 }\,.
\end{eqnarray}
Apart from the universal $R$--matrix ${\cal R}$ that satisfies these
equations for $\Delta_{\cal F}=\Delta ^{op}$ ($\Delta ^{op}=\tau\circ \Delta$, where
$\tau(a\otimes b)=b\otimes a$)  there are two
more well developed  cases of such twists: the jordanian twist of  a Borel algebra
$B(2)$   where ${\cal F}_j$
has the form (\ref{og-twist}) (see \cite{OGIEV}) with $H$ being primitive in
$B(2)$ and $\sigma$ primitive in ${\cal U}_{{\cal F}_j} (B(2))$, and the extended
jordanian  twist (see \cite{KLM} for details).

It will be shown in the next Section that both sets of PE twists are not 
only factorizable but have the factorizable extensions.  One of these 
extensions satisfies ordinary factorization equations (\ref{f-twist1}) and
(\ref{f-twist2}), the other refers to a more sophisticated class.

According to the result by Drinfeld \cite{D3}  skew (constant) 
solutions of the classical Yang--Baxter equation (CYBE) can be quantized and the 
deformed algebras thus  obtained can be presented in a form of twisted universal
enveloping  algebras. On the other hand, such solutions of CYBE can be connected with 
the quasi-Frobenius carrier subalgebras of the initial classical Lie 
algebra \cite{STO}. A Lie  algebra $\g(\mu)$, with the Lie composition $\mu$, is 
called Frobenius if there exists a linear functional $g^* \in \g^*$
such that the form $b(g_1,g_2)=g^*(\mu(g_1,g_2))$ is nondegenerate.
This means that $\g$ must have a nondegenerate 2--coboundary $b(g_1,g_2) \in
B^2(\g,{\bf K})$. The algebra is called quasi-Frobenius if it has a 
nondegenerate 2--cocycle $b(g_1,g_2) \in Z^2(\g,{\bf K})$ (not 
necessarily a coboundary). The classification of quasi-Frobenius
subalgebras in $sl(n)$ can be found in \cite{STO}. In the Section 5
we shall show that extended and peripheric extended twists correspond 
to a class of Frobenius algebras.  

The deformations of quantized algebras include the deformations of their
Lie bialgebras $(\g,\g^*)$. The deformation properties both of $\g$ and 
of $\g^*$ must be taken into consideration. When a Lie algebra 
$\g^*_1(\mu^*_1)$ (with composition $\mu^*_1$) is  deformed  in the 
first order
$$
(\mu^*_1)_t = \mu^*_1 + t\mu^*_2
$$
its deforming function $\mu^*_2$ is also a Lie product and the deformed
property becomes reciprocal: $\mu^*_1$ can be considered as a first order 
deforming function  for algebra $\g^*_2(\mu^*_2)$. Let $\g(\mu)$
be a Lie algebra that form Lie bialgebras both with
$\g^*_1$ and $\g^*_2$. This means that we have a one-dimensional 
family $\{ (\g,(\g^*_1)_t) \}$ of Lie bialgebras and correspondingly a one
dimensional family of quantum deformations $ \{{\cal A}_t(\g,(\g^*_1)_t) \}$ 
\cite{ETI}. This situation provides the possibility to construct in the 
set of Hopf algebras a smooth curve connecting quantizations of the 
type ${\cal A}(\g,\g^*_1)$ with those of ${\cal A}(\g,\g^*_2)$. Such  smooth 
transitions can involve contractions provided
$\mu^*_2 \in B^2 (\g^*_1,\g^*_1)$. This happens in the case of
JT, ET and some other twists (see \cite{KLY} and references therein).


\sect{Extended twists and their limits }

In the construction of extended jordanian twists suggested in \cite{KLM} 
the  carrier algebras of the type ${\bf L}$ play a crucial role. These
are solvable subalgebras with at least four generators.
To study the limit properties of the ET's
let us write down this carrier algebra ${\bf L}$ in the general form:
\begin{equation} \label{el-alg}
\begin{array}{l}
 [H,E] = \delta E, \quad [H,A] = \alpha A, \quad [H,B] = \beta B, \\[0.2cm]
[A,B] = \gamma E, \quad[E,A] = [E,B] = 0,
\end{array}
\end{equation}
\begin{equation} \label{condit}
\alpha + \beta = \delta .
\end{equation}
This parameterization does not describe the full orbit of ${\bf L}$
but presents the essential part of it with the preserved general 
structure of Lie compositions.

In this algebra one can perform successively two nontrivial twists.
The first one corresponds to the carrier subalgebra $B(2)$ with
generators $H$ and $E$. It is called the jordanian twist and has the
twisting element \cite{OGIEV}  
\begin{equation}
\label{ogiev}
\Phi_j = e^{H \otimes \sigma},
\end{equation}
where
\begin{equation}
\label{sigma}
\sigma = \frac{1}{\delta} \ln (1+\gamma E).
\end{equation}
This twisting element is a solution of the factorized twist equations (see
(\ref{f-twist1}) and (\ref{f-twist2})). It transforms the Hopf algebra 
\unal into \unfial 
\begin{equation} 
\label{delfi}
\begin{array}{lcl}
\delfih  & = & H \otimes \exmdes + 1 \otimes H,\\[0.2cm]
\delfia  & = & A \otimes \exas + 1 \otimes A,\\[0.2cm]
\delfib  & = & B \otimes \exbes + 1 \otimes B,\\[0.2cm]
\delfie  & = & E \otimes \exdes + 1 \otimes E.
\end{array}  
\end{equation} 
The jordanian twist (\ref{ogiev}) can be extended \cite{KLM} by the factors
\begin{equation}  
\label{ext-00}
\Phi_{E} = e^{A \otimes B e^{- \beta \sigma }},
\end{equation}
or
\begin{equation}  
\label{ext-01}
\Phi_{E'} = e^{-B \otimes A e^{- \alpha \sigma }}.
\end{equation}
 The element $\Phi_E$ is itself a solution of the general twist equation 
(\ref{gentwist})  for the algebra $U_{j}({\bf L})$. After being twisted 
by $\Phi_E$ the algebra \unfial transforms into \unfal defined by
\begin{equation} 
\label{delf}
\begin{array}{lcl}
\delfh & = & H \otimes \exmdes + 1 \otimes H - \delta A \otimes B 
e^{-(\beta + \delta)\sigma},\\[0.2cm]
\delfa & = & A \otimes \exmbes + 1 \otimes A,\\[0.2cm]
\delfb & = & B \otimes \exbes + \exdes \otimes B,\\[0.2cm]
\delfe & = & E \otimes \exdes + 1 \otimes E.
\end{array}
\end{equation}

The compositions (\ref{el-alg}) and (\ref{delf}) with the condition 
(\ref{condit})
define the three-dimensional set ${\cal H}$ of Hopf algebras. 
All the {\it internal} points of this set correspond to the twisted
algebras of the same general structure and the same properties. To obtain relations
(\ref{delf}) we can also start with \unal  and apply to it the extended twist
${\cal F}_E=\Phi_E\ \Phi_j$ (the composition of $\Phi_E$ and $\Phi_j$). Note also that
for nonzero values of parameters twists
$\Phi_E$ and $\Phi_{E'}$  being applied to algebra
\unalj give the equivalent Hopf algebras $U_{E}({\bf L})\approx U_{E'}({\bf L})$. The
corresponding equivalence map is generated by the substitution $(A,B,\alpha,\beta)
\rightleftharpoons (B,-A,\beta,\alpha)$.

The situation changes when we consider the boundaries of the set ${\cal H}$. 
As we shall see the peripheric Hopf algebras (when they exist) are not only
inequivalent to the initial one but in some cases correspond to a new kind of
extended twists with specific properties.

In the following five cases  the results are trivial: 
\begin{enumerate}
\item $\gamma \rightarrow 0$. The jordanian twist is trivilized. The extensions
become insignificant. They correspond to twisting by primitive elements of an
abelian algebra.   The carrier
subalgebra is here two-dimensional  Abelian and coAbelian.
\item $\delta \rightarrow 0;\, \, \alpha = - \beta \neq 0$. 
In this case the divergences are inevitable in  
\delfa \, and  in \delfb. No limit Hopf algebras in this boundary subset.
\item $ \delta \rightarrow 0$ and $\alpha \rightarrow 0$, $\gamma \neq 0$. 
In such case 
$\beta$ also goes to zero. The behaviour of these parameters can be coordinated
so that the limit Hopf algebra exists (in spite of the divergences of the jordanian
twisting element $\Phi_j$. In this limit the carrier algebra ${\bf L}^{(3)} \equiv
\lim_{\delta,\alpha \rightarrow 0} {\bf L}$   is the central
extension of Heisenberg algebra formed by $A,B$ and $E$. 
 Put $\alpha = a \delta$, $\beta = b\delta$ (with $a + b =1$) and let
 $\sigma_0 \equiv \ln(1 +\gamma E)$. The coproducts of the Hopf algebra
$U_q ({\bf L}^{(3)})$ are defined by the relations
\be
\begin{array}{lcl}
\Delta_q (H) & = & H \otimes e^{-\sigma_0} + 1 \otimes H,\\[1mm]
\Delta_q (A) & = & A \otimes e^{a\sigma_0} + 1 \otimes A,\\[1mm]
\Delta_q (B) & = & B \otimes e^{b\sigma_0} + 1 \otimes B,\\[1mm]
\Delta_q (E) & = & E \otimes e^{\sigma_0} + 1 \otimes E.\\[1mm]
\end{array}
\ee   
Only the last three of them are essential corresponding to 
some special case of Heisenberg algebra quantization. One can easily check
that any group-like elements $f_A, f_B, f^{\prime}_A, f^{\prime}_B$ and $f_E$ depending
on $E$ can serve to construct the coalgebra
\be
\label{heis}
\begin{array}{lcl}
\Delta_q (A) & = & A \otimes f_A +  f^{\prime}_A \otimes A,\\[1mm]
\Delta_q (B) & = & B \otimes  f_B + f^{\prime}_B \otimes B,\\[1mm]
\Delta_q (E) & = & E \otimes f_E + 1 \otimes E,\\[1mm]
\end{array}
\ee
that will form a Hopf algebra with the Heisenberg Lie composition 
$\,[ A,B ] = \gamma E $ in two distinct cases:
\be
\label{first}
f_A f_B = f_E \quad \mbox{\rm and} \quad f^{\prime}_A f^{\prime}_B = 1
\ee
or   
\be
\label{second}
f_A f_B = 1 \quad \mbox{\rm and} \quad f^{\prime}_A f^{\prime}_B = f_E = 1 + 
\widetilde{\gamma}E.
\ee
Thus we have two classes of quantisations of Heisenberg algebra in the scope
of coalgebraic relations (\ref{heis}). The Hopf algebra $U_q ({\bf L}^{(3)})$ 
refers to the first one (with $f^{\prime}_A = f^{\prime}_B =1$ and
$\widetilde{\gamma} = \gamma$). In this case the extensions
\begin{equation}  
\label{ext-02}
\Phi_{E} = e^{A \otimes B f^{-1}_{B}}
\end{equation}
and
\begin{equation}  
\label{ext-03}
\Phi_{E'} = e^{-B \otimes A f^{-1}_{A}}
\end{equation} 
exist and lead to the following quantizations of Heisenberg algebra:
\be
\label{heis-e}
\begin{array}{lcllcl}
\Delta_{q,E} (A) & = & A \otimes f^{-1}_B + 1 \otimes A, 
\quad &\Delta_{q,E^{\prime}} (A) & = & A \otimes f_A + f_{E} \otimes A,\\[1mm]
\Delta_{q,E} (B) & = & B \otimes  f_B + f_E \otimes B,
\quad &\Delta_{q,E^{\prime}} (B) & = & B \otimes f^{-1}_A + 1 \otimes B,\\[1mm]
\Delta_{q,E} (E) & = & E \otimes f_E + 1 \otimes E,
\quad &\Delta_{q,E} (E^{\prime}) & = & E \otimes f_E + 1 \otimes E.\\[1mm]
\end{array}
\ee
Note that $\Delta_{q} (H)$ containing only central elements is not touched
by these extension twists (the same is seen above for $\Delta_{q} (E)$).
Thus the only function of the twists that survive in this case is to bridge 
different classes of quantizations of Heisenberg algebras. 
\item $ \delta \rightarrow 0$ and $\beta \rightarrow 0$. Identical to the previous 
case.
\item $ \delta \rightarrow 0$ and $\gamma \rightarrow 0$.
In this limit the carrier algebra ${\bf L}^{(5)} \equiv \lim_{\delta,\gamma \rightarrow
0} {\bf L}$ is the central extension of the two dimensional algebra $e(2)$ of plane
motions. For the consistent behaviour of parameters the jordanian twist survives
in a form 
$$
\Phi^{(5)}_j = e^{H \otimes \frac{\gamma}{\delta} E }.
$$  
The corresponding deformation $ U({\bf L}^{(5)}) \stackrel{\Phi_j}{\longrightarrow}
U({\bf L}^{(5)}_{j}) $ amounts to a trivial quantization of $U(e(2))$ by a function of
the central generator $E$. No additional transformations are produced by the
extensions $\Phi_E$ or $\Phi_{E^{\prime}}$.

\end {enumerate}
Note that in the second, third and fourth cases the carrier algebra ${\bf L}$ loses the
property of being Frobenius (see Section 5 for more details).

There are two cases that provide nontrivial carrier algebras and twists:

\ i) \underline{ $ \alpha \rightarrow 0;\, \,  \beta = \delta $}.	Let us rewrite the 
corresponding carrier algebra relations:
\begin{equation}
\label{el-alg1}
\begin{array}{l}
\, [H,E] = \delta E, \quad [H,A] = 0, \quad [H,B] = \delta B, \\[2mm]
\, [A,B] = \gamma E, \quad [E,A] = [E,B] = 0.
\end{array}
\end{equation}
This is the limit element of the sequence of algebras of type (\ref{el-alg}), we shall
denote it ${\bf L}^c$. It has rank 2 while  all the other members of the
sequence have rank 1. The twists survive in the limit with the twisting
elements 
\bea
\label{ext-1}
\Phi_j & = & e^{H \otimes \sigma},\\[0.2cm]
\label{ext-11}
\Phi_{P} & = & e^{A \otimes B e^{- \beta \sigma }}.
\eea
The twisted algebra  $U_{j}({\bf L}^c)$ is the
limit of the sequence of Hopf algebras defined by coproducts (\ref{delfi}):
\begin{equation}
\label{delfic1}
\begin{array}{lcl}
\delfih  & = & H \otimes \exmdes + 1 \otimes H,\\[0.2cm]
\delfia  & = & A \otimes 1 + 1 \otimes A,\\[0.2cm]
\delfib  & = & B \otimes \exdes + 1 \otimes B,\\[0.2cm]
\delfie  & = & E \otimes \exdes + 1 \otimes E.
\end{array}  
\end{equation}
The second twisting element $\Phi_{P}$ does not depend on $\delta$ and
leads to the algebra $U_{P}({\bf L}^c)$  with the coproduct:
\begin{equation}
\label{delfc1}
\begin{array}{lcl}
\delfcha & = & H \otimes \exmdes + 1 \otimes H - \delta A \otimes B  
e^{-2\delta \sigma},\\[0.2cm]
\delfcaa & = & A \otimes \exmdes + 1 \otimes A,\\[0.2cm]
\delfcba & = & B \otimes \exdes + \exdes \otimes B,\\[0.2cm]
\delfcea & = & E \otimes \exdes + 1 \otimes E.
\end{array}
\end{equation}
The significant fact is that in $U_{P}({\bf L}^c)$ the element $B\exmdes $ is primitive.
Together with the primitivity of $A$ in 
$U_{j}({\bf L}^c)$ this means that the twisting element $\Phi_{P}$ is now a
solution of the factorized twist equations (\ref{f-twist1}) and (\ref{f-twist2})
contrary to the properties of the internal points of the set
$\widehat{\cal L}$.

ii) \underline{ $ \beta \rightarrow 0;\, \, \alpha = \delta $}. Remember that in the
general situation we have two possible extensions $\Phi_{E}$ and $\Phi_{E'}$ that
give equivalent results. Here the picture is different. On the boundaries of
$\widehat{{\cal L}}$ this degeneracy is removed and we are either to check
both extensions for one type of limits or to study both limits for one of the
extensions. This is the reason to consider this second limit separately.

The purely algebraic part ${\bf L}'^c$ looks like 
\begin{equation}
\label{el-alg2}
\begin{array}{l}
 [H,E] = \delta E, \quad [H,A] = \delta A, \quad [H,B] = 0, \\[0.2cm]
[A,B] = \gamma E, \quad [E,A] = [E,B] = 0 \, ,
\end{array}
\end{equation}
and its jordanian twist $U_{j}({\bf L}'^c)$, 
\begin{equation}
\label{delfic2}
\begin{array}{lcl}
\delfih  & = & H \otimes \exmdes + 1 \otimes H,\\[0.2cm]
\delfia  & = & A \otimes \exdes + 1 \otimes A,\\[0.2cm]
\delfib  & = & B \otimes 1 + 1 \otimes B,\\[0.2cm]
\delfie  & = & E \otimes \exdes + 1 \otimes E,
\end{array}  
\end{equation}
is still equivalent to the previous one, $U_{j}({\bf L}^c)$ (see
(\ref{delfic1})). The extension of the JT has now the form essentially
different from that of  (\ref{ext-11}):
\be
\label{ext-2}
\Phi_{P'}  =  e^{A \otimes B}.
\ee
The final peripheric Hopf algebra $U_{P'}({\bf L}'^c)$ is defined by the
relations:
 \begin{equation}
\label{delfc2}
\begin{array}{lcl}
\delfcc(H) & = & H \otimes \exmdes + 1 \otimes H - \delta A \otimes B  
e^{-\delta \sigma},\\[0.2cm]
\delfcc(A) & = & A \otimes 1 + 1 \otimes A,\\[0.2cm]
\delfcc(B) & = & B \otimes 1 + \exdes \otimes B,\\[0.2cm]
\delfcc(E) & = & E \otimes \exdes + 1 \otimes E.
\end{array}
\end{equation}   
In this case the generator $B$ is primitive in the intermediate algebra (\ref{delfic2})
while $A$ becomes primitive after the extended twist. Thus it does not satisfy the
ordinary factorized twist equations (\ref{f-twist1}) and (\ref{f-twist2}). 
Nevertheless, the relations valid for $\Phi_{P'}$:
\begin{equation}
\label{f-twist3}
\begin{array}{lcl}
(\Delta_{{\cal F}} \otimes {\rm id}) {\cal F} & = & {\cal F}_{13} 
{\cal F}_{23},\\[0.2cm]
 ({\rm id} \otimes \Delta){\cal F} & = & {\cal F}_{12} {\cal F}_{13}.
\end{array}
\end{equation}  
describe the solution of the general twist equation (\ref{gentwist}) in our 
case because both tensor  multipliers in $\Phi_{P'}$ depend each time on 
a single generator providing an  additional commutativity for twisting 
elements in $H \otimes H \otimes H$-space. (Despite the visual similarity
the equations (\ref{f-twist3}) can not be referred to the inverse of the
twisting element ${\cal F}$ due to the structure of the coproduct
$\Delta_{\cal F}$.)

The universal ${\cal R}$-matrices have the form
\be 
{\cal R} =e^{B\exmdes \otimes A }e^{\sigma \otimes H} e^{-H \otimes \sigma}
e^{-A \otimes B \exmdes} 
\ee
in the first case, and
\be 
{\cal R} = e^{B \otimes A}e^{\sigma \otimes H} e^{-H \otimes \sigma} e^{-A \otimes B} 
\ee
in the second.
In both cases the deformation parameter can be introduced by the substitution
$E \rightarrow \xi E; \, A \rightarrow \xi A$. This supplies the deformed
algebra  with the ordinary classical limit  when $\xi \rightarrow 0$, and
gives the  possibility to write down the classical $r$--matrix. It has the
same form in  both cases:
\be
r = A \wedge B + \frac{\gamma}{\delta} H \wedge E,
\ee
(though defined for different carrier algebras (\ref{el-alg1}) and (\ref{el-alg2})). 
Its form guarantees that in both cases the coboundary Lie bialgebras originating 
from it are self-dual.

Just as it was in the case of extended jordanian twist \cite{KLM} one can append any 
number of similar extensions of type $\Phi_{P}$ (correspondingly $\Phi_{P'}$) to the
initial jordanian twist $\Phi_j$ for any number of pairs of equivalent eigenvectors
$(A_m,B_m)$ of the  adjoint operator ${\rm ad}(H)$ and 
with the only nonzero commutators $ [A_m,B_m ] = \gamma E$.

\sect{Peripheric extended twists for simple Lie algebras. $sl(4)$--example.}
 
To demonstrate some other properties of the peripheric extended twists let us apply them
to deform the universal envelopings of simple Lie algebras.
The corresponding carrier subalgebras 
can be found in all the simple Lie algebras with rank no lesser than 2.
We shall work with the algebra $U(sl(4))$ in order to present a completely
nondegenerate case.  The canonical
$gl(4)$--basis $\{ E_{ij} ; \, \,  i,j = 1,\dots,4 \}$ will be used with 
commutation relations
\be
\, [E_{ij},E_{kl}] = \delta_{jk} E_{il} - \delta_{il} E_{kj},
\ee
We shall study the PET with the carrier algebra ${\bf L}'^c$, that is of the
second type (see (\ref{el-alg2})). Let us  injected it into $sl(4)$ in the
following way
\be
\label{el-sl4}
\begin{array}{l}
H = E_{11} - E_{22} \equiv H_{12}, \quad E = E_{24},\\[0.2cm]
A = E_{23}, \qquad\qquad\qquad\ \ \, B = E_{34}.
\end{array}
\ee
This kind of injection corresponds to the fixed values of parameters
\be
\alpha = \delta = -1; \qquad \gamma = 1, \qquad \beta = 0,
\ee  
with
\be
\sigma = -\ln (1 + E_{24}).
\ee 
The universal enveloping algebra $U(sl(4))$ can be twisted by the PET
\be
\label{ext-sl1}
{ {\cal F}_{P'}} = e^{E_{23} \otimes E_{34}} e^{H_{12} \otimes \sigma}.
\ee
The deformed algebra $U_{P'}(sl(4))$ thus obtained has the
comultiplications much less cumbersome compared with the result of an
ordinary ET (see
\cite{KLM}):
\be
\label{twi-co1}
\begin{array}{lcl}
 \delfcc(H_{12}) & = & H_{12} \otimes \exs + E_{23} \otimes E_{34} \exs 
+ 1 \otimes H_{12};\\[0.2cm]
 \delfcc(H_{13}) & = & H_{13} \otimes 1 + 1 \otimes H_{13};\\[0.2cm]
 \delfcc(H_{14}) & = & H_{14} \otimes 1 + 1 \otimes H_{14} + H_{12} \otimes (1 - \exs)
 - E_{23} \otimes E_{34} \exs;\\[0.2cm]
 \delfcc(E_{12}) & = & E_{12} \otimes \ex2s - E_{13} \otimes E_{34} \ex2s 
+ 1 \otimes E_{12} + H_{12} \otimes E_{14}\exs\\ 
		            &   & \   + E_{23} \otimes E_{34} E_{14}\exs;\\[0.2cm]
 \delfcc(E_{13}) & = & E_{13} \otimes \exs + 1 \otimes E_{13} 
- E_{23} \otimes E_{14};\\[0.2cm]
 \delfcc(E_{14}) & = & E_{14} \otimes \exs + 1 \otimes E_{14};\\[0.2cm]
 \delfcc(E_{21}) & = & E_{21} \otimes \exm2s + 1 \otimes E_{21};\\[0.2cm]
 \delfcc(E_{23}) & = & E_{23} \otimes 1 + 1 \otimes E_{23};\\[0.2cm]
 \delfcc(E_{24}) & = & E_{24} \otimes \exms + 1 \otimes E_{24};\\[0.2cm]
 \delfcc(E_{31}) & = & E_{31} \otimes \exms + 1 \otimes E_{31} 
+ E_{21} \otimes E_{34} \exms;\\[0.2cm]
 \delfcc(E_{32}) & = & E_{32} \otimes \exs + 1 \otimes E_{32} 
+ H_{13} \otimes E_{34}\exs;\\[0.2cm]
 \delfcc(E_{34}) & = & E_{34} \otimes 1 + 1 \otimes E_{34} 
+ E_{24} \otimes E_{34};\\[0.2cm]
 \delfcc(E_{41}) & = & E_{41} \otimes \exms + 1 \otimes E_{41} 
+ E_{23} \otimes E_{31}- H_{12} \otimes E_{21}\exs\\ 
		             &   &  \ - E_{23} \otimes E_{34} E_{21} \exs;\\[0.2cm]
\delfcc(E_{42}) & = & E_{42} \otimes \exs - E_{43} \otimes E_{34} \exs + E_{23} \otimes
E_{32}+ 1 \otimes E_{42}\\ 
		             &   & \  - H_{12} \otimes H_{24}\exs +
H_{12} \otimes (\ex2s - \exs)  - E_{23} \otimes H_{24} E_{34} \exs\\
               &   &\  + E_{23} \otimes E_{34}(2\ex2s -\exs) + H_{12}^2 \otimes (\ex2s -
\exs)\\
               &   & \  + 2H_{12}E_{23} \otimes E_{34} \ex2s +
E_{23}^2 \otimes E_{34}^2 \ex2s - H_{12}E_{23} \otimes E_{34} \exs;\\[0.2cm]
\delfcc(E_{43}) & = & E_{43} \otimes 1 + 1 \otimes E_{43} + E_{23} \otimes H_{34}
  - H_{12} \otimes E_{23}\exs\\
               &   &\  - E_{23} \otimes E_{34}E_{23}\exs +
H_{12}E_{23} \otimes E_{24} \exs - E_{23}^2 \otimes E_{34} \exs.
\end{array}
\ee
The  following universal 
${\cal R}$--matrix corresponds to this PET deformation
\be 
{\cal R} = e^{\xi E_{34} \otimes E_{23}}e^{\sigma \otimes H_{12}} e^{-H_{12}
\otimes \sigma}  e^{-\xi E_{23}\otimes E_{34}} .
\ee
In this expression the deformation parameter was introduced (see the previous
section),  so here
$\sigma = -\ln(1+\xi E_{24})$. The corresponding classical $r$--matrix looks
like
\be
r = E_{34} \wedge E_{23} + H_{12} \wedge E_{24}.
\ee

\section{Peripheric twists and Drinfeld--Jimbo quantizations}

It is known for a long time that some types of jordanian quantizations can be treated 
as limit structures for certain smooth sequences of standard deformations
\cite{OGIEV,ABDES,GERST}. It was proved in \cite{KLY} that this property is provided by
the specific correlation between the Lie bialgebras of Drinfeld--Jimbo and ET
quantizations.

Let $(\g, \g^*_{DJ})$ and $(\g, \g^*_{j})$ be the Lie bialgebras corresponding to
Drinfeld--Jimbo and jordanian quantizations of $\g$, respectively. Let $\mu$,
$\mu^*_{DJ}$ and $\mu^*_{j}$  denote the corresponding Lie composition maps. It was
demonstrated in \cite{KLY} that if $\mu^*_{j}$ is a 2--coboundary for the
Lie algebra  $\g^*_{DJ}$, i.e. 
$$\mu^*_{j} \in B^2(
\g^*_{DJ}, \g^*_{DJ}),
$$
then in the set of deformation quantizations of $U(\g)$ there exists a smooth curve
connecting $U_{j}(\g)$ (or in the analogous conditions $U_E(\g)$) with the 
standard deformation $U_{DJ}(\g)$. Smoothness is defined here in the topology very 
similar to the power series one (see \cite{KTOL} and \cite{LYATKA} for details).

It is important to know whether the  algebras twisted by PET's can also be 
connected with DJ quantizations thus describing the limit cases with respect to
the standard deformations. 
In the context of this problem 
we need the inverse of the previous statement. Let us formulate it as follows:

\newtheorem{lemma}{Lemma}
\begin{lemma}\label{lemma1}
Let $U_A(\g)$ and $U_{A'}(\g)$ be two inequivalent quantum deformations of
U(\g) and ${\cal H}(p,q)$ be  a smooth curve  connecting them. If
the curve has the properties:

\ i) ${\cal H}(p,q)_{p=0}= U_{A}(\g)$ and ${\cal H}(p,q)_{q=1}= U_{A'}(\g)$,

ii) ${\cal H}(p,q)$ depends analytically on $q$, 

\noindent then the Lie maps of algebras $\g^*_{A}$ and $\g^*_{A'}$ are
the cocycles for  each other:
\be
\label{nes-cond}
\begin{array}{l}
\mu^*_{A'} \in Z^2(\g^*_{A},\g^*_{A}), \\[0.2cm]
\mu^*_{A} \in Z^2(\g^*_{A'},\g^*_{A'}).
\end{array}
\ee
\end{lemma}
\underline{Proof}. The interior of the set of curves $ \{  {\cal H}(p,q)
, \, q \in \left[ 0,q_1 \right] \}$ forms a neighborhood ${\cal O}(\g)$ of
$U(\g)$ (in the topology induced in the two-dimensional subset ${\cal H}(p,q)$).The
parameters $p$ and $q$ are the natural coordinates in a map covering the 
neighborhood ${\cal O}(\g)$. Thus, for a sufficiently small fixed $q_0  \in \left[ 0,q_1
\right]$ and any small $p$ the pair 
$
\left( \mu \,\, ,q_0\mu^*_{A} + p\mu^*_{A'} \equiv \mu^*_{q_0,p} \right) 
$
is a Lie bialgebra. This means that $\mu^*_{q_0,p}$ is the first order deformation
of $q_0\mu^*_{A}$. But $\mu^*_{A'}$ itself is a Lie algebra. So, $\mu^*_{q_0,p}$
is also the first order deformation of $p\mu^*_{A'}$. $\bullet$

The conditions imposed in this Lemma \ref{lemma1} are natural, they correspond to the 
supposition that there are no singularities in the neighborhood of $U(\g)$ in the
set of its deformation quantizations.  

In the example we have presented in the previous section the Lie map 
$\mu^*_{DJ}(sl(4))$ of the algebra
 $(sl(4))^*_{DJ}$ in the basis $\{ X_{ik} \}$ canonically dual to $\{ E_{ik}
\}$ has the following nonzero commutators:
\be
\label{gdj}
\begin{array}{l}
\, [X_{ii},X_{kl}]_{k \leq l} = \delta_{ik}X_{il} - \delta_{il}X_{ki},\\[0.2cm]
\, [X_{ii},X_{kl}]_{k \geq l} = - \delta_{ik}X_{il} + \delta_{il}X_{ki},\\[0.2cm]
\, [X_{ij},X_{kl}]_{i<j, \, k<l} = 2(\delta_{jk}X_{il} - \delta_{il}X_{kj}),\\[0.2cm]
\, [X_{ij},X_{kl}]_{i>j, \, k>l} = - 2(\delta_{jk}X_{il} - \delta_{il}X_{kj}),
\end{array}
\ee 
The Lie algebra $(sl(4))^*_{P'}$ corresponding to the PET performed by 
(\ref{ext-sl1}) can be extracted from the coproducts (\ref{twi-co1}): 
\be
\label{glet}
\begin{array}{lllll}
\,[X_{11},X_{14}] = & X_{12}, &\quad & [X_{11},X_{21}] = & -X_{41},\\[0.2cm]
\,[X_{11},X_{22}] = & -X_{42}, &  & [X_{11},X_{24}] = & X_{22}-X_{44}, \\[0.2cm]
\,[X_{11},X_{23}] = & -X_{43}, & & [X_{11},X_{34}] = & X_{32}, \\[0.2cm]
\,[X_{11},X_{44}] = & X_{42}, & & [X_{22},X_{21}] = & X_{41},  \\[0.2cm]
\,[X_{22},X_{23}] = & X_{43}, & & [X_{22},X_{24}] = & -X_{22}+X_{44}, \\[0.2cm]
\,[X_{22},X_{44}] = & -X_{42}, & & [X_{33},X_{23}] = & -X_{43}, \\[0.2cm]
\,[X_{33},X_{34}] = & -X_{32},  & & [X_{44},X_{23}] = & X_{43}, \\[0.2cm]
\,[X_{12},X_{24}] = & -2X_{12}, & & [X_{13},X_{24}] = & -X_{13}, \\[0.2cm]
\,[X_{13},X_{34}] = & -X_{12}, & & [X_{14},X_{22}] = & X_{12}, \\[0.2cm]
\,[X_{14},X_{23}] = & X_{13},  & & [X_{14},X_{24}] = & -X_{14}, \\[0.2cm]
\,[X_{21},X_{24}] = & 2X_{21}, & &  [X_{21},X_{34}] = & X_{31}, \\[0.2cm]
\,[X_{23},X_{31}] = & X_{41}, & & [X_{23},X_{32}] = & X_{42},  \\[0.2cm]
\,[X_{23},X_{34}] = & -X_{22}+X_{44}, & &  [X_{24},X_{31}] = & -X_{31}, \\[0.2cm]
\,[X_{24},X_{32}] = & X_{32}, & & [X_{24},X_{34}] = & X_{34}, \\[0.2cm]
\,[X_{24},X_{41}] = & -X_{41},  & & [X_{24},X_{42}] = & X_{42},\\[0.2cm]
\,[X_{34},X_{43}] = & X_{42}. & &  & \\
\end{array}
\ee
We shall denote this set of compositions as $\mu^*_{P'}(sl(4))$. 

One can check by direct computations that the Lie 
multiplications $\mu^*_{DJ}(sl(4))$  and $\mu^*_{P'}(sl(4))$ are not the 
first order deformations of each other. This means (taking into account 
that they are themselves the Lie  compositions) that they are not the  
2--cocycles of each other. So
the conditions (\ref{nes-cond}) are not satisfied and according to the 
Lemma 1 the Hopf algebras $U_{DJ}(sl(4))$ and $U_{P'}(sl(4))$) can not 
be connected by a smooth curve. We have come to the conclusion that
$U_{P'} (sl(4))$) can not be obtained from the Drinfeld--Jimbo 
deformation of $U(sl(4))$ by a contraction or by any other smooth 
limit process. This feature clearly shows how different could be the 
results of quantum deformations by extended and by peripheric twists.

The facts discussed above are tightly connected with 
the problem of equivalence of different CYBE solutions and in this
context with the properties of the corresponding quasi-Frobenius
algebras. We have seen that all the algebras belonging to the set 
$\widetilde{\cal L}=\{{\bf L}(\alpha,\delta - \alpha,\gamma,\delta)|\gamma \neq 0,
\delta \neq 0 \}$ are at least quasi-Frobenius. This property can be
precised as follows.
\begin{lemma}\label{lemma2}
All the elements of the set $\widetilde{\cal L}$ are Frobenius algebras.  
\end{lemma}
\underline{Proof}. For all the algebras {\bf L} of the set $\widetilde{\cal
L}$ the form  
$$
b(g_1,g_2) = E^*(\left[ g_1,g_2 \right]) \quad \quad g_1, g_2 \in {\bf L}
$$
is nondegenerate. Here $E^*$ is the functional 
canonically dual to the basic element $E \in {\bf L}$. 

\ \hfill$\bullet$

Note that our results are in total agreement with the classification of
quasi-Frobenius algebras of low dimensions given by Stolin \cite{STO}.
One can check that the set $\widetilde{\cal L}$ is equivalent to the 
class $\{  P_{a_1,a_2,a_3} | a_1 \neq a_3 \}$ 
(see the Proposition 1.2.3 in \cite{STO}). 

\section{Conclusions}

The peripheric twists described in this paper are not continuously
connected with Drinfeld--Jimbo deformations despite the fact that the
carrier subalgebras of the peripheric and ordinary extended twists
belong to the same smooth family of Frobenius algebras. Taking into account
that the $U_{E}(sl(n))$ algebra quantized by certain types of ET can be 
treated as continuous limit of DJ deformations \cite{KLY} we have at least 
the superposition of two smooth transitions that can connect DJ and PET
deformations. In the case studied above the algebra ${\bf
L}(\alpha,0,\gamma,\delta) \subset sl(4)$ can be obtained from  ${\bf
L}(1,1,1,2) \subset sl(3) \subset sl(4)$ by means of  a ``rotation" in the 
space of the Cartan subalgebra of $sl(4)$. We want to stress that the
``rotation" connecting ${\bf L}(1,1,1,2)$ with ${\bf L}(-1,0,1,-1)$ is not a
similarity transformation for ${\bf L}$ and thus cannot be used to carry
properties from ET to PET and viceversa. Nevertheless it might be possible to
simulate analogous "rotations" also in the set of modified DJ deformations
(using multiparametric quantizations or applying the continuous families
of dual groups \cite{LYMU}). If both "rotations" could be matched there
might exist the possibility of contraction-like smooth transition between      
modified DJ and PET deformations.   

\section*{Aknowlegements}   

The authors are thankful to Prof. P.P.Kulish for his important comments.
V.L. would like to thank the DGICYT of the Ministerio de Educaci\'on y Cultura de 
Espa\~na   for supporting his sabbatical stay (grant SAB1995-0610). This work has been
partially supported by DGES of the  Ministerio de  Educaci\'on y Cultura of Espa\~na
under Project PB95-0719, the Junta de Castilla y Le\'on (Espa\~na) and the Russian 
Foundation for Fundamental Research under grant  97-01-01152.


\end{document}